\documentclass[fleqn]{mat01}
\usepackage{times,mathtimy,amssymb,latexsym}
\begin{document}

\setcounter{page}{477}
\firstpage{477}

\def\d{\mbox{\rm d}}

\def\C{\cal C}
\def\A{\cal A}
\def\sm{\setminus}
\def\ro{\rho}
\def\bigint{\int}
\def\ss{\subset}
\def\B{\cal B}
\def\D{\cal D}
\def\U{\cal U}
\def\V{\cal V}
\def\R{\Bbb R}
\def\N{\Bbb N}
\def\T{\Bbb T}
\def\L{\Bbb L}
\def\phi{\varphi}
\def\e{\varepsilon}
\def\H{\Bbb H}
\def\P{\cal P}
\def\O{\Bbb O}
\def\Q{\Bbb Q}
\def\G{\Bbb G}
\def\exp{\operatorname{exp}}
\def\id{\operatorname{id}}
\def\bs{\setminus}
\def\Comp{{\cal Comp}}
\def\M{\cal M}
\def\Auth{\operatorname{Auth}}
\def\min{\operatorname{min}}
\def\max{\operatorname{max}}
\def\inf{\operatorname{inf}}
\def\diam{\operatorname{diam}}
\def\Cl{\operatorname{Cl}}
\def\Fr{\operatorname{Fr}}
\def\Int{\operatorname{Int}}

\newtheorem{theore}{Theorem}
\renewcommand\thetheore{\arabic{section}.\arabic{theore}}
\newtheorem{theor}[theore]{\bf Theorem}
\newtheorem{lem}[theore]{Lemma}
\newtheorem{coro}[theore]{\rm COROLLARY}

\newtheorem{theoree}{Theorem}
\renewcommand\thetheoree{\arabic{theore}.\arabic{theoree}}
\newtheorem{poft}[theoree]{\it Proof of Theorem}

\title{On topological properties of the Hartman--Mycielski functor}

\markboth{Taras Radul and Du\v{s}an Repov\v{s}}{On topological
properties of the Hartman--Mycielski functor}

\author{TARAS RADUL and DU\v{S}AN REPOV\v{S}$^{*}$}

\address{Departmento de Matematicas, Facultad de Cs. Fisicas y Matematicas,
Universidad de Concepcion, Casilla 160-C, Concepcion, Chile\\
\noindent $^{*}$Institute for Mathematics, Physics and Mechanics,
University of Ljubljana, Jadranska 19, Ljubljana, Slovenia 1001\\
\noindent E-mail: tarasradul@yahoo.co.uk; dusan.repovs@fmf.uni-lj.si}

\volume{115}

\mon{November}

\parts{4}

\pubyear{2005}

\Date{MS received 2 January 2005}

\begin{abstract}
We investigate some topological properties of a normal functor $H$
introduced earlier by Radul which is some functorial
compactification of the Hartman--Mycielski construction HM. We
prove that the pair ($HX$, HM$Y$) is homeomorphic to the pair
$(Q,\sigma)$ for each nondegenerated metrizable compactum $X$ and
each dense $\sigma$-compact subset $Y$.
\end{abstract}


\maketitle

\section{Introduction}

The general theory of functors acting on the category $\Comp$ of
compact Hausdorff spaces (compacta) and continuous mappings was
founded by Shchepin \cite{Sh}. He described some elementary
properties of such functors and defined the notion of the normal
functor which has become very fruitful. The classes of all normal
and weakly normal functors include many classical constructions:
the hyperspace exp, the space of probability measures $P$, the
superextension $\lambda$, the space of hyperspaces of inclusion
$G$, and many other functors (see \cite{FZ} and \cite{TZ}).

Let $X$ be a space and $d$ an admissible metric on $X$ bounded by
$1$. By $\hbox{HM}X$ we shall denote the space of all maps from
$[0,1)$ to the space~$X$ such that $f|[t_i,t_{i+1})\equiv{}$
const, for some $0=t_0\leq\cdots\leq t_n=1$, with respect to the
following metric:
\begin{equation*}
d_{\rm HM}(f,g) = \bigint_0^1 d(f(t),g(t))\d t,\quad f,g\in
\hbox{HM}X.
\end{equation*}

The construction of $\hbox{HM}X$ is known as the
Hartman--Mycielski construction \cite{HM}. Recall, that the
Hilbert cube is denoted by $Q$, and the following subspace of $Q$
\begin{equation*}
\{(a_n)_{n=1}^\infty\in Q\mid a_k=0 \ \ \text{for all but finitely
many }k\}
\end{equation*}
is denoted by $\sigma$. Telejko has shown in \cite{Te} that for
any nondegenerated separable metrizable $\sigma$-compact strongly
countable-dimensional space $X$ the space HM$X$ is homeomorphic to
$\sigma$.

For every $Z\in\Comp$ consider
\begin{align*}
\hbox{HM}_nZ = \big\{f\in \hbox{HM}Z\mid &\hbox{ there exist
} 0=t_1<\cdots<t_{n+1}=1\\[.2pc]
&\hbox{ with }f|[t_i,t_{i+1})\equiv z_i\in Z, i=1,\dots,n\big\}.
\end{align*}

Let $\U$ be the unique uniformity of $Z$. For every $U\in\U$ and
$\e>0$, let
\begin{equation*}
\langle \alpha,U,\e \rangle = \{\beta\in \hbox{HM}_nZ\mid
m\{t\in [0,1)\mid (\alpha(t),\beta(t'))\notin U\}<\e\}.
\end{equation*}
The sets $\langle \alpha,U,\e\rangle$ form a base of a compact
Hausdorff topology in $\hbox{HM}_n Z$. Given a map $f:X\rightarrow Y$ in
$\Comp$, define a map $\hbox{HM}_n X\rightarrow \hbox{HM}_n Y$ by the
formula $\hbox{HM}_n F(\alpha)=f\circ\alpha$. Then $\hbox{HM}_n$
is a normal functor in $\Comp$ (see \cite{TZ}, \S2.5.2).

For $X\in\Comp$ we consider the space $\hbox{HM}X$ with the
topology described above. In general, $\hbox{HM}X$ is not compact.
Zarichnyi has asked if there exists a normal functor in $\Comp$
which contains all functors $\hbox{HM}_n$ as subfunctors
\cite{TZ}. Such a functor $H$ was constructed in \cite{Ra}.

We investigate some topological properties of the space $HX$ which
is some natural compactification of $\hbox{HM}X$. The main results
of this paper are as follows:

\begin{theor}[\!]
$HX$ is homeomorphic to the Hilbert cube for each nondegenerated
metrizable compactum $X$.
\end{theor}

\begin{theor}[\!]
The pair $(HX, {\rm HM}Y)$ is homeomorphic to the pair
$(Q,\sigma)$ for each nondegenerated metrizable compactum $X$ and
each dense $\sigma$-compact strongly countable-dimensional subset
$Y$.
\end{theor}

\section{Construction of \pmb{$H$}}

Let $X\in\Comp$. By $C(X)$ we denote the Banach space of all
continuous functions $\phi:X\rightarrow\R$ with the usual $\sup$-norm:
$\|\phi\| =\sup\{|\phi(x)|\mid x\in X\}$. We denote the segment
$[0,1]$ by $I$.

For $X\in\Comp$ let us define the uniformity of $\hbox{HM}X$. For
each $\phi\in C(X)$ and $a,b\in [0,1]$ with $a<b$ we define the
function $\phi_{(a,b)}\!:\hbox{HM}X\rightarrow \R$ by
\begin{equation*}
\phi_{(a,b)}=\frac 1{(b-a)}\bigint_a^b\phi\circ\alpha(t)\d t.
\end{equation*}
Define
\begin{equation*}
S_{\rm HM}(X)=\{\phi_{(a,b)}\mid \phi\in C(X) \quad \hbox{and}
\quad (a,b)\subset [0,1)\}.
\end{equation*}

For $\phi_1,\dots,\phi_n\in S_{\rm HM}(X)$ define a pseudometric
$\ro_{\phi_1,\dots,\phi_n}$ on $\hbox{HM}X$ by the formula
\begin{equation*}
\ro_{\phi_1,\dots,\phi_n}(f,g)=\max\{|\phi_i(f)-\phi_i(g)|\mid
i\in\{1,\dots,n\}\},
\end{equation*}
where $f,g\in \hbox{HM}X$. The family of pseudometrics
\begin{equation*}
\P=\{\ro_{\phi_1,\dots,\phi_n}\mid n\in\N, \quad \hbox{where} \
\phi_1,\dots,\phi_n\in S_{\rm HM}(X)\},
\end{equation*}
defines a totally bounded uniformity $\U_{{\rm HM}X}$ of
$\hbox{HM}X$ \cite{Ra}.

For each compactum $X$ we consider the uniform space
$(HX,\U_{HX})$ which is the completion of $(\hbox{HM}X,\U_{{\rm
HM}X})$ and the topological space $HX$ with the topology induced
by the uniformity $\U_{HX}$. Since $\U_{{\rm HM}X}$ is totally
bounded, the space $HX$ is compact.

Let $f\!\!:X\rightarrow Y$ be a continuous map. Define the map $\hbox{HM} f:
\hbox{HM}X\rightarrow \hbox{HM}Y$ by the formula
$\hbox{HM}f(\alpha)=f\circ\alpha$, for all $\alpha\in \hbox{HM}X$.
It was shown in \cite{Ra} that the map
$\hbox{HM}f:(\hbox{HM}X,\U_{{\rm HM}X})\rightarrow(\hbox{HM}Y,\U_{{\rm
HM}Y})$ is uniformly continuous. Hence there exists the continuous
map $Hf:HX\rightarrow HY$ such that $Hf|\hbox{HM}X=\hbox{HM}f$. It is easy
to see that $H:\Comp\rightarrow\Comp$ is a covariant functor and
$\hbox{HM}_n$ is a subfunctor of $H$ for each $n\in\N$.

\section{Preliminaries}

All spaces are assumed to be metrizable and separable. We begin
this section with the investigation of certain structures of
equiconnectivity on $HX$ for some compactum $X$. The first one
will be the usual convexity.

Let us remark that the family of functions $S_{\rm HM}(X)$ embed
$\hbox{HM}X$ in the product of closed intervals
$\prod_{\phi_{(a,b)}\in S_{\rm HM}(X)}I_{\phi_{(a,b)}}$, where
$I_{\phi_{(a,b)}}=[\min_{x\in X} |\phi(x)|,\max_{x\in X}$
$|\phi(x)|]$. Thus, the space $HX$ is the closure of the image of
$\hbox{HM}X$. We denote by $p_{\phi_{(a,b)}}:HX\rightarrow
I_{\phi_{(a,b)}}$ the restriction of the natural projection.

\setcounter{theore}{0}
\begin{lem}
$HX$ is a convex subset of $\prod_{\phi_{(a,b)}\in S_{\rm
HM}(X)}I_{\phi_{(a,b)}}$.
\end{lem}

\begin{proof}
It is enough to prove that for each $\alpha$, $\beta\in
\hbox{HM}X\subset \prod_{\phi_{(a,b)}\in S_{\rm
HM}(X)}I_{\phi_{(a,b)}}$ we have that
$\frac12\alpha+\frac12\beta\in HX$.

Let $0=r_0<r_1<\cdots<r_k=1$ and $0=p_0<p_1<\cdots<p_m=1$ be the
decompositions of the unit interval which correspond to $\alpha$
and $\beta$. Consider any $\phi^1,\dots,\phi^n\in C(X)$ and
$(a_i,b_i)\subset (0,1)$ for $i\in\{1,\dots,n\}$. Choose a
decomposition of the unit interval $0=t_0<t_1<\cdots<t_s=1$ such
that $\{r_0,\dots,r_k\}\cup\{p_0,\dots,p_m\}\cup\{a_i|i
\in\{1,\dots,n\}\}\cup\{b_i|i\in\{1,\dots,n\}\}\subset
\{t_0,\dots,t_s\}$.

Define $\gamma\in \hbox{HM}X$ as follows. Let $t\in[t_l,t_{l+1})$
for some $l\in\{0,\dots,s-1\}$. Define $\gamma(t)=\alpha(t)$ if
$t<\frac12(t_l+t_{l+1})$ and $\beta(t)$ otherwise. It is easy to
see that $p_{\phi^i_{(a_i,b_i)}}(\gamma)=
p_{\phi^i_{(a_i,b_i)}}(\frac12\alpha+\frac12\beta)$ for each
$i\in\{1,\dots,n\}$. Hence $\frac12\alpha+\frac12\beta\in
\Cl(\hbox{HM}X)= HX$. The lemma is proved.
\end{proof}

\begin{coro}$\left.\right.$\vspace{.5pc}

\noindent $HX$ is absolute retract for each metrizable compactum
$X$.
\end{coro}

Define the map $e_1:\hbox{HM}X\times \hbox{HM}X\times I\rightarrow
\hbox{HM}X$ by the condition that $e_1(\alpha_1,\alpha_2,t)(l)$ is
equal to $\alpha_1(l)$ if $l<t$ and $\alpha_2(l)$ in the opposite
case for $\alpha_1,\alpha_2\in \hbox{HM}X$, $t\in I$ and $l\in
[0,1)$. We consider $\hbox{HM}X$ with the uniformity $\U_{{\rm
HM}X}$ and $I$ with the natural\break metric.

\begin{lem}
The map $e_1:\hbox{HM}X\times \hbox{HM}X\times I\rightarrow \hbox{HM}X$ is
uniformly continuous.
\end{lem}

\begin{proof}
Let us consider any $U\in\U_{{\rm HM}X}$. We can suppose that
\begin{equation*}
U=\{(\alpha,\beta)\in \hbox{HM}X\times \hbox{HM}X\mid
|\phi_{(0,1)}(\alpha)-\phi_{(0,1)}(\beta)|<\delta\},
\end{equation*}
for some $\delta>0$ and $\phi\in C(X)$. The proof of the general
case is the same.

Put $c=\max_{x\in X} |\phi(x)|$. Choose $n\in\N$ such that
$1/n<\delta/2c$ and put $a_i=i/n$ for $i\in\{0,\dots,n\}$.
Consider an element $V$ of uniformity $\U_{{\rm HM}X}$ defined as
follows:
\begin{align*}
V &= \left\lbrace (\alpha,\beta)\in \hbox{HM}X\times \hbox{HM}X \mid
|\phi_{(a_i,a_{i+1})}(\alpha)-\phi_{(a_i,a_{i+1})}(\beta)|<
\frac{\delta}{2n^2},\right.\\
&\qquad \qquad \hbox{for each} \ \ i\in\{0,\dots,n-1\}\bigg\}.
\end{align*}
Put $E=\{(l,s)\in I\times
I\mid|l-s|<\frac{1}{2n}\}$.

Let us consider any $((\alpha_1,\beta_1),(\alpha_2,\beta_2),
(l,s))\in V\times V\times E$. Then we have
\begin{align*}
&\left|\int_{a_i}^{a_{i+1}}(\phi\circ\alpha_1(t)-\phi\circ\alpha_2(t))\d t
\right|\\[.3pc]
&\quad\, = n\frac{1}{a_{i+1}-a_i} \left| \int_{a_i}^{a_{i+1}}
(\phi\circ\alpha_1(t)-\phi\circ\alpha_2(t))\d t \right|\\[.3pc]
&\quad\, = n|\phi_{(a_i,a_{i+1})}(\alpha_1)-\phi_{(a_i,a_{i+1})}
(\alpha_2)|<n\frac{\delta}{2n^2}=\frac{\delta}{2n},
\end{align*}
for each $i\in\{0,\dots,n-1\}$. We have the same for $\beta_1$
and $\beta_2$. Since $|t_1-t_2|<\frac{1}{2n}$, there exists
$i_0\in\{0,\dots,n-1\}$ such that $t_1,t_2\in
[a_{i_0},a_{i_0+1}]$. Then we have
\begin{align*}
&|\phi_{(0,1)}(e_1(\alpha_1,\beta_1,t_1))-\phi_{(0,1)}(e_1(\alpha_2,\beta_2,t_2))|\\[.3pc]
&\quad\, = \left|\int_0^1(\phi\circ e_1(\alpha_1,\beta_1,t_1)(t) -
\phi\circ e_1(\alpha_2,\beta_2,t_2)(t))\d t \right|\\[.3pc]
&\quad\, \leq \sum_{i=0}^{i_0-1}
\left|\int_{a_i}^{a_{i+1}}(\phi\circ
\alpha_1(t) - \phi\circ \alpha_2(t))\d t \right| + c(a_{i_0+1}-a_{i_0})\\[.3pc]
&\qquad\, +\sum_{i=i_0}^{n-1} \left|\int_{a_i}^{a_{i+1}}(\phi\circ
\beta_1(t)-\phi\circ \beta_2(t))\d t \right|
<\frac{n-1}{2n}\delta+\frac{\delta}{2}<\delta.
\end{align*}
Therefore $(e_1(\alpha_1,\beta_1,t_1), e_1(\alpha_2,
\beta_2,t_2))\in V$. Thus the lemma is proved.
\end{proof}

Hence there exists the extension of $e_1$ to the continuous map
$e:HX\times HX\times I\rightarrow HX$. It is easy to check that
$e(\alpha,\alpha,t)=\alpha$, for each $\alpha\in HX$.

\section{Proofs}

\setcounter{theore}{1}
\begin{poft}{\rm
Since $HX$ is infinite-dimensional convex compactum, Theorem~1.1
follows by Klee Theorem \cite{Kl}.}
\end{poft}

We will use the equiconnectivity defined by the function $e$ to
prove Theorem~1.2. Let $(Z,e)$ be an equiconnected space where
$e:Z\times Z\times I\rightarrow Z$ is the map which defines the structure
of equiconnectedness. A subset $A\subset Z$ is called $e$-{\it
convex} if $e(a,b,t)\in A$, for each $a$, $b\in A$ and $t\in I$.

Consider any $n\in\N$. Let $D_{n-1}$ be a standard $n-1$-simplex
in $\R^n$. Define a map $e_n:Z^n\times D_{n-1}\rightarrow Z$ as follows.
Consider any $(a_1,\dots,a_n,\lambda_1,\dots,\lambda_n)\in
Z^n\times D_{n-1}$. Define the finite sequence
$\{x_1,\dots,x_n\}$ by induction. Put $x_1=a_1$. Suppose that we
have already defined $x_j$ for each $j\leq i-1$ where $1<i\leq n$.
Define
\begin{equation*}
x_i=e \left(x_{i-1},a_i, \frac{\sum_{l=1}^{i-1}
\lambda_l}{\sum_{l=1}^i\lambda_l} \right).
\end{equation*}
Put $e_n(a_1,\dots,a_n,\lambda_1,\dots,\lambda_n)=x_n$. One can
check that the function $e_n$ is continuous.

Let us recall that $A\subset Z$ is homotopically dense in $Z$ if
there exists a homotopy $H:Z\times I\rightarrow Z$ such that $H(Z\times
(0,1])\subset A$.

\setcounter{theore}{0}
\begin{lem}
Let $(Z,e)$ be an equiconnected absolute retract and let $A\subset
Z$ be a dense $e$-convex subset. Then A is homotopically dense in
$Z$.
\end{lem}

\begin{proof}
It is enough to prove that for every $n\in\N$, every point $z\in
Z$ and every neighborhood $U\subset Z$ of $z$ there exists a
neighborhood $V$ of $z$ such that every map $f:\partial I^n\rightarrow
V\cap A$ admits an extension $F:I^n\rightarrow U\cap A$ \cite{To}.

Since $e_{n+1}(z,\dots,z,\lambda_1,\dots,\lambda_{n+1})=z$ for
each $(\lambda_1,\dots,\lambda_{n+1})\in D_n$, there exists a
neighborhood $V$ of $z$ such that
$e_{n+1}(z_1,\dots,z_{n+1},\lambda_1,\dots,\lambda_{n+1})\in U$
for each $z_1,\dots,z_{n+1}\in V$ and
$(\lambda_1,\dots,\lambda_{n+1})\in D_n$.

Consider any map $f:\partial I^n\rightarrow V\cap A$. Choose any metric
$d$ on $I^n$. Let us consider a Dugunji system for $I^n\sm\partial
I^n$, that is, an indexed family $\{U_s,a_s\}_{s\in S}$ such that

\begin{enumerate}
\renewcommand\labelenumi{(\arabic{enumi})}
\item $U_s\subset I^n\sm\partial I^n$, $a_s\in \partial I^n$ ($s\in
S$),

\item $\{U_s\}_{s\in S}$ is a locally finite cover of
$I^n\sm\partial I^n$,

\item If $x\in U_s$, then $d(x,a_s)\leq 2d(x,\partial I^n)$ for
$s\in S$ (see of Ch.~II, \S3 of \cite{BP}) for more details.
\end{enumerate}

Since $\dim I^n=n$, we can suppose that for each $x\in
I^n\sm\partial I^n$ there exists a neighborhood of $X$ which meets
at most $n+1$ elements of $\{U_s\}_{s\in S}$. Moreover we can
suppose that the index set is countable. Choose some partition of
unity $\{b_s\}_{s\in S}$ inscribed into $\{U_s\}_{s\in S}$. Fix
some order $S=\{s_1,s_2,\dots\}$. Define a function $F:I^n\rightarrow U$
as follows. Put $F(x)=f(x)$ for each $x\in\partial I^n$. Now
consider any $x\in I^n\sm\partial I^n$. There exists a finite
sequence $(s_{i_1},\dots,s_{i_{n+1}})$ such that $x\notin U_s$
for each $s\in S\sm\{s_{i_1},\dots,s_{i_{n+1}}\}$. Put $F(x)=
e_{n+1}(f(a_{s_{i_1}}),\dots,f(a_{s_{i_{n+1}}}),b_{s_{i_1}}(x),\dots,
b_{s_{i_{n+1}}}(x))$. One can check that the function $F$ is a
continuous extension of the function $f$. Since $A$ is $e$-convex,
it follows that $f(I^n)\subset A$. The lemma is thus proved.
\end{proof}

\begin{poft}{\rm
It is easy to check that $\hbox{HM}X$ is an $e$-convex subset of
$HX$. Since $HX$ is an equiconnected absolute retract, and
$\hbox{HM}Y$ is a dense $e$-convex subset, homeomorphic to
$\sigma$ \cite{Te}, Theorem~1.2 follows by \cite{BRZ}; Proposition~3.1.7 and
Lemma~4.1.}
\end{poft}

\section*{Acknowledgement}

This research was supported by the research grant SLO-UKR 02-03/04.

\end{document}